\documentclass[11pt, reqno]{amsart}
\usepackage{amsmath, amsfonts, amsthm, amssymb, multicol, mathtools, dsfont,verbatim}
\usepackage{graphicx}
\usepackage{float, hyperref}
\usepackage{scalerel,stackengine, subcaption}
\usepackage[usenames,dvipsnames,x11names]{xcolor}
\usepackage{enumitem}
\usepackage{pgfplots}
\usepgfplotslibrary{fillbetween}
\pgfplotsset{width=10cm,compat=1.9}
\usepackage[square,sort,comma,numbers]{natbib}
\makeatletter
\def\@setauthors{%
  \begingroup
  \def\thanks{\protect\thanks@warning}%
  \trivlist
  \centering\footnotesize \@topsep30\p@\relax
  \advance\@topsep by -\baselineskip
  \item\relax
  \author@andify\authors
  \def\\{\protect\linebreak}

  \normalsize\lowercase{\authors}%
  
	\ifx\@empty\contribs
  \else
    ,\penalty-3 \space \@setcontribs
    \@closetoccontribs
  \fi
  \endtrivlist
  \endgroup
}
\def\@settitle{\begin{center}
\LARGE\lowercase{\@title}
  \end{center}%
}
\makeatother
\newcommand{\authoremail}[1]{\email{\href{mailto:#1}{\color{lightblue}{#1}}}}
\newcommand{\authoraddress}[1]{\address{\normalfont{#1}}}

\usepackage{fancyhdr}
\numberwithin{equation}{section}

\theoremstyle{definition}

\theoremstyle{remark}

\renewcommand{\epsilon}{\varepsilon}

\renewcommand{\geq}{\geqslant}
\renewcommand{\leq}{\leqslant}

\newcommand{\ubd}{\dim_{\textup{B}}}

\newcommand{\ad}{\dim_{\mathrm{A}} }

\newcommand{\as}{\dim^\theta_{\mathrm{A}} }

\newcommand{\hd}{\dim_{\mathrm{H}}  }

\newcommand{\be}{\begin{equation}}
\newcommand{\ee}{\end{equation}}

\renewcommand{\epsilon}{\varepsilon}
\newcommand{\eps}{\varepsilon}

\newcommand{\rd}{\mathbb{R}^d}

\newcommand{\fs}{\dim^\theta_{\mathrm{F}}}

\makeatletter
\DeclareRobustCommand\widecheck[1]{{\mathpalette\@widecheck{#1}}}
\def\@widecheck#1#2{%
    \setbox\z@\hbox{\m@th$#1#2$}%
    \setbox\tw@\hbox{\m@th$#1%
       \widehat{%
          \vrule\@width\z@\@height\ht\z@
          \vrule\@height\z@\@width\wd\z@}$}%
    \dp\tw@-\ht\z@
    \@tempdima\ht\z@ \advance\@tempdima2\ht\tw@ \divide\@tempdima\thr@@
    \setbox\tw@\hbox{%
       \raise\@tempdima\hbox{\scalebox{1}[-1]{\lower\@tempdima\box
\tw@}}}%
    {\ooalign{\box\tw@ \cr \box\z@}}}
\makeatother

\stackMath
\newcommand\reallywidehat[1]{%
\savestack{\tmpbox}{\stretchto{%
  \scaleto{%
    \scalerel*[\widthof{\ensuremath{#1}}]{\kern.1pt\mathchar"0362\kern.1pt}%
    {\rule{0ex}{\textheight}}
  }{\textheight}%
}{2.4ex}}%
\stackon[-6.9pt]{#1}{\tmpbox}%
}
\definecolor{lightblue}{HTML}{2B77A4}
\colorlet{plotblue}{LightSkyBlue3!80}
\definecolor{darkred}{HTML}{9E0D0D}
\definecolor{darkyellow}{HTML}{b3b300}
\definecolor{darkorange}{HTML}{D86129}
\hypersetup{
	colorlinks=true,
	linkcolor=darkred,
	urlcolor=darkred,
	citecolor=lightblue
}
\title{An invitation to dimension interpolation}

\author{Jonathan M. Fraser}
\authoraddress{J. M. Fraser, University of St Andrews, Scotland}
\authoremail{jmf32@st-andrews.ac.uk}
\thanks{JMF was financially supported by an \emph{EPSRC Open Fellowship} (EP/Z533440/1), \emph{A new theory of dimension interpolation} and  a  \emph{Leverhulme Trust Research Project Grant} (RPG-2023-281).}

\begin{document}
\maketitle
\thispagestyle{empty}

\begin{abstract}
A \emph{fractal} is an object   exhibiting complexity at arbitrarily small scales.  In order to study and characterise fractals, one is often interested in quantifying how they fill up space on small scales. This gives rise to various notions of \emph{fractal dimension}.  However, even for the simplest examples, the different  definitions of dimension may completely disagree about the answer.  In this expository article I will examine this phenomenon   and  use it to discuss and motivate \emph{dimension interpolation}. 
Dimension interpolation views these classical notions as boundary points of continuous families of dimensions, thus transforming isolated numerical answers into a coherent geometric picture.
\end{abstract}

\tableofcontents

\section{Dimension theory} \label{intro}

Dimension is one of the most fundamental concepts in mathematics.  A line is 1-dimensional, a filled in square is 2-dimensional, a solid cube is 3-dimensional, and so on.  Fractals, however,  complicate this picture and often refuse to conform to our  intuition about dimension. They  are objects exhibiting complexity at arbitrarily small scales and appear across mathematics and wider science.  Their somewhat vague definition  precludes smooth objects such as circles or smooth curves because they eventually  look like a line upon magnification---they have a `tangent'.  Fractals, on the other hand, will keep changing and producing more complexity under magnification.
\[
\textit{How should we define the dimension of a fractal?}
\]
Defining `dimension' in a theoretically robust and geometrically meaningful way is already a nuanced problem---and one of my favourite problems  in mathematics. There  are many distinct notions of `fractal dimension' and many (but not all!) of these involve quantifying how the fractal object fills up space on small scales.  Here we present three natural and distinct ways of doing this and in doing so an apparent paradox will emerge:~even for the simplest possible example, our definitions of dimension will give completely different answers!   

In order to understand how a fractal $X \subseteq \rd$ fills up space on small scales, we first  fix a scale, let's say $r>0$ which we think of as being very small (later we will let it tend to zero).  Fixing a scale serves to  discretise the problem and allows us to start to calculate.
\[
\textit{How big is   $X$ at scale $r$?  }
\]
In order to answer this question, we count how many balls of diameter $r$ are required to contain (or `cover') $X$.  That is, how small can I make $n$ such that there are balls $B_1, \dots, B_n$ (or arbitrary  sets) of diameter $r$ such that
\[
X \subseteq \bigcup_{k=1}^n B_k.
\]
Such a collection $\{B_k\}_k$ is known as an \emph{$r$-cover} of $X$ or simply a \emph{cover}.   Denote the smallest such $n$ by $N_r(X)$ to emphasise the dependence on $r$ and $X$. The  number $N_r(X)$ then captures how big $X$ is at scale $r$.
\[
\textit{How should we expect $N_r(X)$ to behave?  }
\]
We are interested in the behaviour of $N_r(X)$ as $r \to 0$ since this will produce a scale invariant description of how $X$ fills up space on small scales. It is important to insist that our notions of dimension are `scale invariant' because, for example,  a line segment is 1-dimensional no matter how long it is.  What distinguishes the size of different line segments is their length, but length is an inherently 1-dimensional notion of size, making dimension a more fundamental concept than length.

 It turns out to be quite difficult (and ultimately not very important) to compute  $N_r(X)$ exactly.  But what we can often derive is the \emph{growth rate} of $N_r(X)$ as $r \to 0$, that is, establish expressions of the form
\begin{equation} \label{covering}
N_r(X) \lesssim  r^{-s}
\end{equation}
for some exponent $s>0$.  The use of $\lesssim$ in \eqref{covering} means that there exists a    constant $C$ (perhaps depending on $X$ and $s$ but \emph{not} on $r$) such that $N_r(X) \leq C  r^{-s}$.  This is excellent notation for this kind of problem (and we will use it throughout this article) because it suppresses unnecessary details and emphasises that it is the $s$ we care about, not the $C$.  The growth estimate \eqref{covering} becomes harder to satisfy as $s$ becomes smaller and the smallest $s$ we can get away with can be thought of as the `dimension' of $X$.  For example
\[
N_r([0,1]) \lesssim  r^{-1}, \qquad N_r([0,1]^2) \lesssim r^{-2}, \qquad N_r([0,1]^3) \lesssim r^{-3}
\]
and these exponents cannot be improved.  This is consistent with our idea that a line is 1-dimensional, a square is 2-dimensional, and a cube is 3-dimensional.  Formally, the \emph{(upper) box dimension} $\ubd X$ of a bounded set $X$ is the infimum of $s \geq 0$ such that for some $C>0$ and all $r>0$, 
\[
N_r(X) \leq C r^{-s}.
\]
The box dimension is geometrically meaningful, elegantly defined, and clearly describes what we were trying to capture. 

 Problem solved? Well, maybe not quite.  There are theoretical drawbacks of this definition.  For example, since we assume the balls in the cover are the same size and we naively count the number required, we do not account for the possibility that the covers used are wildly inefficient.  That is, we may find a much more efficient cover by using balls  with \emph{different} diameters.  However,  if we do this we cannot simply count the number of covering sets, and must instead weight the contribution of each covering  set by its size and then consider the total `cost' of a cover in a more nuanced way.  This approach gives rise to the Hausdorff dimension.  Formally, the \emph{Hausdorff dimension}  $\hd X$ of $X$ is defined as the infimum of exponents $s \geq 0$ such that for all $\eps>0$ there exists  a cover $\{B_k\}_k$ (by balls with possibly \emph{different} diameters denoted by $|B_k|$) such that
 \[
 \sum_k |B_k|^s < \eps.
 \]
 That is,  the smallest $s$ such that $X$ appears small when viewed from  an $s$-dimensional perspective.  For $s>\ubd X$, taking an efficient $r$-cover by balls of the \emph{same} size yields
 \[
 \sum_{k=1}^n |B_k|^s = N_r(X) \,  r^s \to 0
 \]
 and so 
 \[
 \hd X \leq \ubd X
 \]
 always holds, but we shall see that the inequality can be strict! The Hausdorff dimension is  theoretically robust, and satisfies many natural and elegant properties which a good notion of dimension should.
 
 Problem solved? Well, maybe not quite.  Both the Hausdorff and box dimensions give rise to \emph{global} notions of dimension which provide a description of how $X$ fills up space \emph{on average}.  While very useful, this does not paint the whole picture.  Often this information is \emph{not} what is needed for a particular application, which may be sensitive to \emph{local} or \emph{extreme}  information. This consideration leads to the Assouad dimension.  First, we localise at scale $R>0$ and location $x$ and then we apply the `box-counting approach' at scale $r<R$ \emph{but only} to $X \cap B(x,R)$, where $B(x,R)$ is the ball centred at $x$ of radius $R$.  The localised covering number now depends on the `localisation scale' $R>0$,  location $x$,  and the `covering scale' $r<R$, and must be interpreted with this in mind.  Formally, the \emph{Assouad dimension} $\ad X$ of $X$ is defined as the infimum of exponents $s \geq 0$ such that there exists $C>0$ such that for all $R>r>0$ and $x$
 \[
 N_r\Big( X \cap B(x,R) \Big) \leq C \left(\frac{R}{r} \right)^s.
 \]
 That is,  the smallest $s$ such that  $X$ looks small \emph{at all scales and locations} when viewed from an $s$-dimensional perspective.  For bounded $X$, setting $R$ equal to  the diameter of $X$, we immediately see that
 \[
 \ubd X \leq \ad X
 \]
  always holds, but we shall again see that the inequality can be strict! The Assouad dimension gives clear and readily applicable geometric information.
 
So, are we finished? No! There are many further nuances and many other notions of dimension which we will not discuss here, but more relevant to our discussion is that we now have three perfectly natural, geometrically meaningful, theoretically relevant notions of fractal dimension, all of which describe how $X$ fills up space on small scales.  
 \[
 \emph{How do we decide which dimension to choose?}
 \]
   The answer to this question is not to choose!  We should consider them \emph{all} and interpret them carefully and in the context of what the definition was trying to capture.   For more background on box and Hausdorff dimensions see \cite{Fal03} and  for more on the Assouad dimension see \cite{jon:book, robinson}.
 
 \section{A worked example}

In this section we  consider a simple but very instructive example.  In fact, the example will be the simplest object which one could possibly think of as a fractal.  Let 
\[
X = \big\{  1, \, 1/2,\,  1/3, \, 1/4, \, 1/5, \,  \dots\big\} = \{ 1/n : n \in \mathbb{N}\} .
\]
That is, $X$ is the set of reciprocals of  positive integers.  It is a countable discrete set with a single accumulation point at 0 and everything is very explicit.  What could be easier?
 
We consider each notion of dimension  in turn, starting with the smallest. First, the Hausdorff dimension of $X$  is as small as possible, that is, 
\[
\hd X = 0.
\]
This can be seen by taking the cover $\{B_k\}_k$ given by
\[
B_k = B\big(1/k, 2^{-k-1} \delta \big)
\]
for some small $\delta>0$, that is, $B_k$ is an interval of diameter $2^{-k} \delta$ centred at the point $1/k \in X$.  Then, for all positive $s>0$, 
 \[
 \sum_k |B_k|^s = \sum_k 2^{-ks}\delta^s = \frac{2^{-s}\delta^s}{1-2^{-s}}  
 \]
which can be made arbitrarily small by choosing $\delta$ sufficiently small.  In fact this argument works for \emph{any} countable set $X$, showing that all countable sets have Hausdorff dimension zero.
 
 Next we show that the box dimension is different from the Hausdorff dimension and, in fact,
 \[
 \ubd X = 1/2.
 \]
 Given $0<r<1/2$, choose $k \in \mathbb{N}$ (uniquely) such that 
 \[
 \frac{1}{k} - \frac{1}{k+1} \leq r <  \frac{1}{k-1} - \frac{1}{k} .
 \]
 This separates the points in $X$ into two distinct classes.  The first $k$ points are separated by more than $r$ from each other which means that each point will require its own covering set.  Therefore,
  \[
 N_r(X) \geq   k \gtrsim r^{-1/2}
 \]
 and so the box dimension is \emph{at least 1/2}.   On the other hand, naively covering each of the first $k$ points separately and then covering the rest crudely gives
 \begin{equation} \label{boxdimest}
 N_r(X) \leq k+N_r(X \cap [0,1/k])  \lesssim k+ \frac{1/k}{r} \lesssim r^{-1/2},
 \end{equation}
 which shows that the box dimension is \emph{at most} 1/2 and  proves the claim. 
 
 Finally, consider the Assouad dimension.  Let $x=0$ and $r=R^2$ with $R$ small. Then all  the points in $X$ lying within $R$ of $x=0$ are within $r=R^2$ of their neighbours, and therefore
 \[
 N_r(X \cap B(0,R)) \gtrsim N_r([0,R]) \approx \frac{R}{r}
 \]
 and so the Assouad dimension is \emph{at least} 1 but in fact
 \[
 \ad X = 1
 \]
 since $X$ lives in the line and the line itself clearly has Assouad dimension 1. Indeed, we can always naively apply the bound
 \[
  N_r(X \cap B(x,R)) \lesssim N_r([x-R,x+R]) \lesssim \frac{R}{r}
 \]
 to obtain an upper bound of 1.
 
 We are now in quite  a remarkable situation.  Even for the \emph{simplest possible fractal}, the three notions of dimension considered here \emph{completely disagree} with each other.  Not only do they return three distinct candidates for the dimension, but these candidates are as uniformly spread and separated from each other as they could possibly be, given that $X$ is a subset of the line.   
\[
\textit{How should we interpret this striking disagreement?}
\]
The key point is that these are three \emph{different perspectives} and each gives a clear geometric interpretation of $X$.  The Hausdorff dimension thinks that $X$ is as small as possible since, by considering highly inhomogeneous covers, the set can be covered with negligible cost.  The box dimension, on the other hand, thinks that $X$ has an intermediate size, because it recognises that, at scale $r$, $X$ coarsely fills up a significant amount of space, while still recognising the substantial `holes' appearing in the set far from the origin.  Finally, the Assouad dimension thinks $X$ is as large as possible because when one localises around the accumulation point, $X$ maximally fills up space at the scale given by the square of the localisation scale:~the line itself would be no larger!

None of the answers given by the distinct notions of dimension are wrong; rather, each makes sense from its own perspective. An appreciation of this sort of `disagreement' goes a long way, even outside of mathematics.

\section{Dimension interpolation}

\emph{Dimension interpolation} is a relatively novel programme in fractal geometry, which aims to address and interpret the type of situation we  observed in the previous section. It suggests a new holistic way to think about dimension:~rather than treat the different notions of dimension in isolation, we should attempt to view them as different facets of a single object.  The idea is to define \emph{interpolating families} of dimensions which deform continuously from one notion to another.  If this can be done, then we will be able to witness one perspective changing into another, and we will be able to reconcile the differences of opinion expressed by the distinct notions of dimension.  Moreover, we will uncover richer and more nuanced geometric information about the fractal set in question and be able to leverage this information in various applications, both specific to fractal geometry and in wider mathematics.  

Here we discuss two examples of dimension interpolation in action. The first example is the \emph{intermediate dimensions} which were introduced by Falconer, Fraser and Kempton in \cite{intdims} to interpolate between the Hausdorff   and box dimensions. With the interpolation goal in mind,  first observe  that the   difference between the definitions of Hausdorff and box dimension is that for Hausdorff dimension there is \emph{no restriction} on the relative sizes of  diameters of sets used in the cover, but for the box dimension all covering sets must have the \emph{same diameter}.    The intermediate dimensions are then defined by using an interpolation parameter $\theta \in [0,1]$ to restrict the range of allowable diameters used in the covers.  It turns out that the most interesting way to do this is to insist that,  for any two sets $U$ and $V$ in the cover, $|U| \leq |V|^\theta$, that is, all covering sets have diameters in the range $[\delta^{1/\theta}, \delta]$ for some $\delta \in (0,1)$.

More precisely, the  \emph{(upper) intermediate dimension} $\dim_\theta X$ is defined as the infimum of exponents $s \geq 0$ such that for all $\eps>0$ there exists  a cover $\{B_k\}_k$   such that
 \[
 \sum_k |B_k|^s < \eps
 \]
 and such that for all $k,l$,  $|B_k| \leq |B_l|^\theta$.
 That is,  the smallest $s$ such that $X$ appears small when viewed from  an $s$-dimensional and $\theta$-restricted perspective.    If $\theta=0$, then there are no restrictions on the diameters and the definition recovers the Hausdorff dimension.  On the other hand, if $\theta=1$, then all covering sets are forced to have the same diameter and 
the cost reduces to
 \[
 \sum_k |B_k|^s  = N_r(X) \, r^s
 \]
 and the definition recovers the box dimension. 
 
 Various properties of intermediate dimensions are established in \cite{intdims}. In particular, $\dim_\theta X$ is monotonically increasing in $\theta\in [0,1]$, is  continuous except perhaps at $\theta = 0$,  is invariant under bi-Lipschitz mappings, and satisfies
 \[
 \hd X \leq \dim_\theta X \leq \ubd X
 \]
 for all $\theta$.  The possible functions which can be realised as the intermediate dimensions of some set were completely classified in \cite{intclass}.  It turns out that there is a very rich family of possibilities!  
 
 The intermediate dimensions have already found numerous applications.  These applications include the dimension theory of Brownian images \cite{burrell}, bi-Lipschitz classification problems and   multifractal analysis \cite{banajicarpets},  Sobolev mapping problems \cite{frasertyson}, and  the dimension theory of orthogonal projections   \cite{bff21}.

The second example of dimension interpolation in action is the \emph{Assouad spectrum} which was  introduced by   Fraser and Yu in \cite{assouadspectrum} to interpolate between box  dimension and Assouad dimension. This was the first appearance of `dimension interpolation'.  Again with the interpolation goal in mind, first observe  that the definition of  Assouad dimension uses  \emph{two scales} (the `localisation scale' $R$ and the `covering scale' $r$) but the box dimension just uses \emph{one scale} (the `covering scale' $r$). If the covering scale is extremely small relative to the localisation scale, then covering a small piece is not so different from covering the whole set, and one expects the  box dimension to emerge.  Moreover, pairs of scales which witness the Assouad dimension are often closer together than one might expect (with the obvious caveat that $R/r \to \infty$).  The Assouad spectrum is defined by using an  interpolation parameter $\theta \in (0,1)$  to fix the relationship between these two scales.  It turns out that the most interesting way to do this is to set $R = r^\theta$. 

More precisely, the \emph{Assouad spectrum} $\as X$ of $X$ is defined as the infimum of exponents $s \geq 0$ such that there exists $C>0$ such that for all $0<r<1$ and $x$
 \[
 N_r\left( X \cap B\big(x,r^\theta\big) \right) \leq C \left(\frac{r^\theta}{r} \right)^s.
 \]
 That is, the smallest $s$ such that $X$ looks small \emph{at all scales and locations} when viewed from an $s$-dimensional and $\theta$-restricted perspective.  If $\theta = 0$ and $X$ is bounded, then one immediately gets $\dim_{\textup{A}}^0X = \ubd X$ and as $\theta \to 1$ the localisation and covering scales become relatively close to each other and this is how the Assouad dimension is typically witnessed. In many cases (but not all!) we get
 \[
 \as X \to \ad X 
 \]
as $\theta \to 1$. Various properties of the Assouad spectrum were established in \cite{assouadspectrum} including that it is continuous in $\theta \in [0,1)$, is invariant under bi-Lipschitz maps, and satisfies
\[
  \ubd X \leq  \as X \leq     \ad X
\]
for $\theta \in (0,1)$. The possible functions which can be realised as the Assouad spectrum  of some set were completely classified in \cite{specclass}.  Once again, it turns out that there is a very rich family of possibilities!

 The Assouad spectrum has already found numerous applications, some coming in surprising areas.  These applications include $L^p \to L^q$ mapping properties of certain operators from  harmonic analysis \cite{roos, beltran}, weak embeddability problems  \cite{stathis},    quasiconformal mapping problems \cite{stathistyson}, and  the Sullivan dictionary from conformal dynamics \cite{bullams}.

\section{Our worked example:~revisited}

Recall the simple example described above, which drew highly disparate opinions from our three notions of dimension. That is, the prototypical fractal,
\[
X =\{1/n\}_{n \in \mathbb{N}}.
\]
Of course, the reader is now wondering how dimension interpolation handles this example. We might expect a continuous interpolation through the three distinct notions of dimension and that the form of this interpolation   reveals geometric information not available to the individual notions---indeed, this will be the case!  One can show, with a little more work this time, that for all $\theta \in (0,1)$
\[
 \dim_\theta X = \frac{ \theta}{1+\theta} \qquad  \quad \text{and} \qquad  \quad  \as X = \min\left\{ \frac{1/2}{1-\theta}, \,  1\right\}.
\]
Thus the intermediate dimensions increase smoothly from 0 to $1/2$, while the Assouad spectrum increases from $1/2$ to 1, with a phase transition at $\theta = 1/2$. A dimension theory sceptic might look at this simple example and declare that it is countable and therefore the Hausdorff dimension is zero, it is not a fractal, and there is nothing more to say.   But this is the wrong perspective:~even in this very simple case, an appreciation of the different perspectives on dimension and dimension interpolation has unearthed much more geometric information about $X$.  Moreover, we uncover a complete interpolation between 0 and 1, which observes  the three isolated notions of fractal dimension in the process.  We are also led to many further questions pertaining to finer geometric features of $X$:~why are the intermediate dimensions strictly concave?  why does the Assouad spectrum have a phase transition at $\theta=1/2$? why does the Assouad spectrum reach the Assouad dimension for values of $\theta<1$?  Are these behaviours typical for other examples? and so on.

\begin{figure}[H] 
	\centering
	\includegraphics[width=0.9\textwidth]{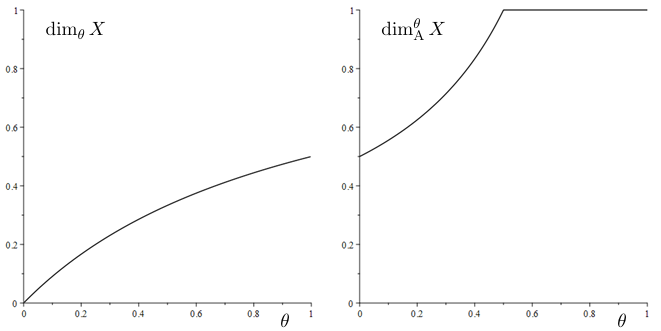}
\caption{\emph{Interpolation in action}:~plots of the   intermediate dimensions (left), and the Assouad spectrum (right) as functions of $\theta \in (0,1)$ for the simple example $X =\{1/n\}_{n \in \mathbb{N}}$. }
\end{figure}

We conclude this article by proving the formulae presented above for the intermediate dimensions and Assouad spectrum of our simple example. The proofs are straightforward and self-contained but will require some subtle observations. More general results  are available in \cite{intdims,assouadspectrum}, respectively.

Let's begin with  the intermediate dimensions.  We first bound  $\dim_\theta X$ above, which we do by constructing explicit covers with diameters permissible according to the $\theta$-restriction in the definition.  In fact we only use two distinct scales $r$ and $r^\theta$.  Let $s>0$ be arbitrary for now, $0<r<1$ and let $M =\lceil r^{-(s +\theta(1-s))/2}\rceil$.   Take a cover $\{B_k\}_k$ of  $X$ consisting of $M$ intervals $B_k=B(k^{-1}, r/2)$  of length $r$ for $1\leq k\leq M$ and $\lesssim M^{-1}/r^\theta $ intervals  of length $r^\theta$  that cover $[0, M^{-1}].$   Then
\begin{eqnarray} 
\sum_{k} |B_k|^s &\lesssim & Mr ^s +  \frac{r^{\theta s}}{Mr^\theta} \label{balance}\\ 
 & \lesssim &  r^{-(s +\theta(1-s))/2}  r^s 
+ r^{\theta (s-1)}r^{(s +\theta(1-s))/2} \nonumber \\
& =&  2r^{(s(\theta+1)-\theta)/2} \ \to \ 0 \nonumber
\end{eqnarray}
as $r \to 0$ provided  $s(\theta +1) > \theta$. Thus $\dim_\theta X \leq\theta/(1+\theta)$.  In case the reader is wondering where the choice of $M$ came from, one can first use an arbitrary $M$ and then seek to balance the two terms on the right hand side of \eqref{balance}.

For the lower bound we show $X$ is large by showing we may spread mass out thinly on $X$ in an appropriate sense.  This is an intermediate dimension analogue of the well-known \emph{mass distribution principle} used to study the Hausdorff dimension, see \cite{Fal03}.

Fix $s = \theta/(1+\theta)$. Let $0<r < 1$ and again let $M =\lceil r^{-(s +\theta(1-s))/2}\rceil$, which now simplifies as $M =\lceil r^{-s}\rceil$.  Define a mass distribution  $\mu_r$ on $X$ as the sum of point masses   with
\begin{equation*}\label{mass}
 \mu_r \big(\left\{1/k\right\}\big)\  = \ 
\left\{
\begin{array}{cl}
 r^s & \mbox{ if } 1\leq k \leq  M   \\
 0 &    \mbox{ if } M+1\leq k <\infty  
\end{array}
\right. .
\end{equation*}
Then
\begin{eqnarray*} 
\mu_r(X) =Mr^s \geq 1 .
\end{eqnarray*}
 We now want to show that $\mu_r$ is `thinly spread out' by showing that $\mu_r (B)$ is small for sets  $B$  such that $r \leq |B|\leq r^\theta$. To this end, note that for $k \leq M$,
$$\frac{1}{(k-1)} - \frac{1}{k}\ \geq \ \frac{1}{k^{2}}\ \geq \ \frac{1}{M^{2}},$$
and thus the gap between any two points of $X$ carrying mass is at least $1/M^{2}$.
Then such sets  $B$ intersect at most 
$1+ |B|/(1/M^{2})  = 1+|B|M^{2}  \lesssim |B| r^{-2s}$ of the points of $X$ which have mass $r^s$.
Hence 
\begin{eqnarray*} 
\mu_r (B) \lesssim   r^s|B| r^{-2s} \leq |B|^s
\end{eqnarray*}
where we used that $|B| \leq r^\theta = r^{s/(1-s)}$. This is our mass distribution estimate. Now, let $\{B_k\}_k$ be a cover of $X$ such that $r \leq |B_k|\leq r^\theta$ for all $k$. Then
$$1\ \leq \ \mu_r(X) \ \leq\  \mu_r\Big(\bigcup_k B_k\Big)\ \leq \ \sum_k \mu_r(B_k) \
\lesssim \   \sum_k |B_k|^s,$$
and so $\sum_k |B_k|^s \gtrsim 1$ for every admissible cover and therefore $\dim_\theta X \geq s$.

Finally, consider the Assouad spectrum. By naively bounding the covering number of a small ball by the covering number of the whole space and then applying the box dimension estimate \eqref{boxdimest}, for $0<r<1$ and $x$, we obtain
\[
 N_r\left( X \cap B\big(x,r^\theta\big) \right) \leq N_r(X) \lesssim r^{-1/2} = \left(\frac{r^\theta}{r}\right)^{\frac{1/2}{1-\theta}}
\]
which proves
\[
\as X \leq \frac{1/2}{1-\theta}.
\]
The bound $\as X \leq 1$ is easy since $X$ lives in the line and so 
\[
 N_r\left( X \cap B\big(x,r^\theta \big) \right)  \lesssim  \frac{r^\theta}{r}.
\]
For the lower bound, we already proved that
\[
\dim_{\textup{A}}^{1/2} X = 1
\]
because when we proved that $\ad X = 1$ above we witnessed this using scales $r=R^2$, that is, we used $\theta=1/2$.  It is then straightforward to see that $\as X \geq 1$ for all $\theta \in [1/2,1)$.  For $\theta \in (0,1/2)$, let $r \in (0,1/2)$, $x=0$, and choose $k\in \mathbb{N}$  uniquely  such that 
 \[
 \frac{1}{k} - \frac{1}{k+1} \leq r <  \frac{1}{k-1} - \frac{1}{k} 
\]
noting that $k \approx r^{-1/2}$ and so $1/k < r^\theta$.  The points in $X \cap [0,1/k] \subseteq B(0,r^\theta)$ are all closer than $r$ to their neighbours and therefore covering $X \cap [0,1/k]$ at scale $r$ is no easier than covering $[0,1/k]$ at scale $r$.  Therefore
\[
 N_r\left( X \cap B\big(0,r^\theta\big) \right) \gtrsim   \frac{1/k}{r} \approx  r^{-1/2} = \left(\frac{r^\theta}{r}\right)^{\frac{1/2}{1-\theta}}
\]
which proves 
\[
\as X \geq \frac{1/2}{1-\theta}
\]
as required.

\section{More interpolation}

We refer those interested in more detailed treatments of dimension interpolation to the survey articles \cite{Fra19, intprojsurvey}.  There are of course many other ways to interpolate between notions of fractal dimension---and many different questions one can ask about the notions we have discussed in this short article.  For example, the \emph{Fourier spectrum} was introduced in \cite{fourierspectrum} and interpolates between the Fourier and Hausdorff dimensions.  For more background on Fourier dimension, and its many connections with Hausdorff dimension, see \cite{mattilafourier}.  The Fourier spectrum $\fs X$ is non-decreasing, continuous for $\theta \in [0,1]$, and recovers the Fourier and Hausdorff dimensions precisely at the end points.  Moreover, it has been used to study the well-known \emph{Fourier restriction problem} in harmonic analysis \cite{restriction} and the \emph{Falconer distance problem} from geometric measure theory  \cite{fourierspectrum}.

Dimension interpolation is still in its infancy.  However, the emerging richness of the theory and the---often surprising---applications  suggest that it is not merely a refinement of existing notions, but a genuinely new lens through which geometric  phenomena can be studied. Looking into the future, we may start to see the  classical notions of dimension as merely boundary points of richer families.  Moreover, the  additional geometric information uncovered by these families may be crucial to applications in an increasing  variety of areas.


\end{document}